\documentclass[12pt,twoside]{article}
\usepackage{amsmath,amsthm,amssymb}
\usepackage{amsfonts}
\usepackage{mathrsfs}
\usepackage{graphicx}
\usepackage[a4paper]{geometry}
\begin{document}
\title{\textbf{Some results on the norm of finite groups}}
\author{Mark L. Lewis, Zhencai Shen and Quanfu Yan}
\date{}

\maketitle

\begin{abstract}
Let $G$ be a finite group and $N_{\Omega}(G)$ be the intersection of the normalizers of all subgroups belonging to the set $\Omega(G),$ where $\Omega(G)$ is a set of all subgroups of $G$ which have some theoretical group property. In this paper, we show that $N_{\Omega}(G)= Z_{\infty}(G)$ if $\Omega(G)$ is one of the following: (i) the set of all self-normalizing subgroups of $G$; (ii)  the set of all subgroups of $G$ satisfying the subnormalizer condition in $G$; (iii) the set of all pronormal subgroups of $G$; (iv)  the set of all $\mathscr{H}$-subgroups of $G$; (v) the set of all  weakly normal subgroups of $G$; (vi)  the set of all $NE$-subgroups of $G$.
\end{abstract} 

\section{Introduction}

All the groups considered in this paper are finite and $G$ denotes a finite group.The notation and terminology used in this paper are
standard, as in \cite{Huppert}
In 1935, R. Baer \cite{Baer1} considered the norm $N(G)$ of $G$, which is the intersection of normalizers of all subgroups of $G$.
In 1953, R. Baer \cite{Baer2} considered the intersection of normalizers of all Sylow subgroups of $G$ and proved that it coincides with the hypercenter $Z_{\infty}(G)$, where $Z_{\infty}(G)$ is the terminal term of the upper central series of the group $G$. \\

Let $\Omega(G)$ be the set of all subgroups of the group $G$ which have some theoretical group property. Define $N_{\Omega}(G)$, called the $\Omega$-norm $G$, be the intersection of the normalizers of all subgroups belonging to the set $\Omega(G).$ In the case $\Omega(G)=\emptyset,$ we assume that $G=N_{\Omega}(G).$  It is clear that the norm $N(G)$ is the $\Omega$-norm of $G$ for the set $\Omega(G)$, which consists of all subgroups of $G$, and $Z_{\infty}(G)$ is the $\Omega$-norm of $G$ for the set $\Omega(G)$, which consists of all Sylow subgroups of $G$.  \\

In this paper, we prove the following: \\

\noindent{\bf Theorem 1.1.} Let $G$ be a group. Then $N_{\mathcal{SN}}(G)=Z_{\infty}(G),$ where $\mathcal{SN}(G)$ is the set of all self-normalizing subgroups of $G$.\\

Notice that the norm $N_{\mathcal{SN}}(G)$ can be viewed as the intersection of all self-normalizing subgroups of $G.$ Similarly, we consider the intersection of all self-centralizing subgroups of $G$ and get the following:\\

\noindent{\bf Theorem 1.2.} Let $G$ be a group. Then $Z(G)$ coincides with the intersection, denoted by $C_{\mathcal{SC}}(G)$, of all subgroups in $\mathcal{SC}(G)$, where $\mathcal{SC}(G)$ is the set of all self-centralizing subgroups of $G$.\\

\noindent{\bf Remark 1.3.} Notice that both $N_{\mathcal{SC}}(G)=Z_{\infty}(G)$ and $N_{\mathcal{SC}}(G)=Z(G)$ are not always true. For example, let $G=Q_{16}$ be the generalised quaternion group of order 16. Then $Z(G)<N_{\mathcal{SC}}(G)=Z_{2}(G)<Z_{3}(G)=G.$ This is the reason why we consider the intersection $C_{\mathcal{SC}}(G)$ of all self-centralizing subgroups of $G$.\\ 

\noindent{\bf Theorem 1.4.} Let $G$ be a group and $p$ be a prime. Let $\Omega_p(G)$ be the set of all $p$-groups in $\Omega(G),$ where $\Omega(G)$ is one of the following:\\
(i) the set of all pronormal subgroups of $G$.\\ 
(ii) the set of all $\mathscr{H}$-subgroups of $G$.\\
(iii) the set of all  weakly normal subgroups of $G$.\\
(iv)  the set of all subgroups of $G$ satisfying the subnormalizer condition in $G$.\\
(v) the set of all $NE$-subgroups of $G$.\\
Then\\
\noindent(a) $N_{\Omega}(G)$ coincides with the intersection of normalizers of all sylow subgroups of $G.$\\
(b) if $\Omega_p(G)$ is not the case (v), $N_{\Omega_p}(G)$ coincides with the intersection of normalizers of all sylow $p$-subgroups of $G$.\\

\section{Definitions}

In this section, we recall some definitions and introduce some notation.\\

 A subgroup $H$ of a group $G$ is called pronormal in $G$ if, for any $g \in G$, there exists $x \in \langle H, H^g \rangle$ such that $H^x = H^g.$ Bianchi in \cite{Bianchi} introduced the concept of $\mathscr{H}$-subgroups as follows:\\
 
\noindent{\bf Definition 2.1.} A subgroup $H$ of a group $G$ is an $\mathscr{H}$-subgroup of $G$ if, for any $g \in G$, $N_G(H)\cap H^g \leq H.$ \\

In \cite{Muller}, M\"uller introduces the weak normality: \\

\noindent{\bf Definition 2.2.} A subgroup $H$ of $G$ is called weakly normal in $G$ if $H^g \leq N_G(H)$  implies that $g \in N_G(H).$ \\

In \cite{Shirong}, Shirong Li introduced the definition of $NE$-subgroups of $G$.\\ 

\noindent{\bf Definition 2.3.} A subgroup $H$ is called an $NE$-subgroup of $G$ if $H=N_G(H)\cap H^G,$ where $H^G$ is the normal closure of $H$ in $G.$\\

In \cite{Mysovskikh}, Mysovskikh introduced the following subgroup embedding property.\\ 

\noindent{\bf Definition 2.4.} A subgroup $H$ of $G$ is said to satisfy the subnormalizer condition in $G$ if for every subgroup $K$ of $G$ such that $H \unlhd K$, it follows that $N_G(K) \leq N_G(H).$ \\

For a prime $p$, we define the following sets of subgroups of $G$:\\

\noindent (i) $\mathcal{S}(G)$ is the set of all Sylow subgroups of $G$.\\
(ii) $\mathcal{SN}(G)$ is the set of all self-normalizing subgroups of $G$, that is, $\mathcal{SF}(G)=\{H\leq G| H=N_G(H)\}.$ \\
(iii) $\mathcal{SC}(G)$ is the set of all self-centralizing subgroups of $G$, that is, $\mathcal{SC}(G)=\{H\leq G| C_G(H)\leq H\}.$ \\
(iv) $\Sigma(G)$ is the set of all subgroups of $G$ satisfying the subnormalizer condition in $G$.\\
(v) $\mathcal{S}_p(G)$ is the set of all Sylow $p$-subgroups of $G$.\\
(vi) $\Sigma_p(G)$ is the set of all $p$-subgroups of $G$ which belong to $\Sigma(G)$. \\


\section{Proofs of Main Results}

To prove Theorem 1.1, we begin by the following obvious lemma.\\

\noindent{\bf Lemma 3.1.} Let $G$ be a group. Assume that $N$ is normal in $G$ and $N\leq H\leq G.$ Then $H\in\mathcal{SN}(G)$ if and only if $H/N\in\mathcal{SN}(G/N).$\\

\noindent{\bf Proof of Theorem 1.3.} Notice that $Z_{\infty}(G)=\cap_{P\in\mathcal{S}(G)}N_G(P)$ by Baer's result and $N_G(P)=N_G(N_G(P))\in\mathcal{SN}(G).$ Then $N_{\mathcal{SN}}(G)\leq Z_{\infty}(G).$ Thus, we only need to show that $Z_{\infty}(G)\leq N_{\mathcal{SN}}(G).$\\

We work by induction on the order of $G.$ If $Z(G)=1,$ then $Z_{\infty}(G)=1$ and so $N_{\mathcal{SN}}(G)=1$, we are done. Assume that $Z(G)>1.$ Then $Z_{\infty}(G/Z(G))\leq N_{\mathcal{SN}}(G/Z(G))$ by induction. Notice that if $H\in\mathcal{SN}(G),$ then $Z(G)\leq N_G(H)=H$. Hence $H/Z(G)\in\mathcal{SN}(G/Z(G))$ if and only if $H\in\mathcal{SN}(G)$ by Lemma 3.1. It follows that 
\begin{eqnarray*}    
N_{\mathcal{SN}}(G/Z(G))&=&\bigcap_{H/Z(G)\in\mathcal{SN}(G/Z(G))}H/Z(G)  \nonumber    \\
&=&  \bigcap_{Z(G)\in H\in\mathcal{SN}(G)}H/Z(G)\\
&=&\bigcap_{H\in\mathcal{SN}(G)}H/Z(G)\\
&=&(\bigcap_{H\in\mathcal{SN}}H)/Z(G)\\
&=& N_{\mathcal{SN}}(G)/Z(G).
\end{eqnarray*}

On the other hand, we have $Z_{\infty}(G)/Z(G)=Z_{\infty}(G/Z(G)).$ It follows that $Z_{\infty}(G)/Z(G)\leq N_{\mathcal{SN}}(G)/Z(G)$ and so $Z_{\infty}(G)\leq N_{\mathcal{SN}}(G),$ as wanted. $\Box$ \\

The proof of Theorem 1.2 requires the following lemma.\\

\noindent{\bf Lemma 3.2.} Let $A$ be a maximal abelian subgroup of a group $G$.  Then $A=C_G(A)$ and so $A$ is self-centralizing in $G$. \\

\noindent{\bf Proof.} It is clear that $A\leq C_G(A).$ Assume that $A<C_G(A).$ Then there is an element $x\in C_G(A)$ but $x\not\in A.$ Consider the group $\langle A, x \rangle$ generated by $A$ and $x$. As $x\in C_G(A),$ we have $\langle A, x \rangle$ is abelian, contrary to the maximal choice of $A$. Hence $A=C_G(A)$, which is self-centralizing by definition. $\Box$ \\

\noindent{\bf Proof of Theorem 1.2.}  Let $H$ be a self-centralizing subgroup of $G$. Notice that $Z(G)\leq C_G(H)\leq H$ and so $Z(G)\leq C_{\mathcal{SC}}(G).$\\
Let $g\in C_{\mathcal{SC}}(G).$ Let $A$ be any maximal abelian subgroups of $G$. Then by Lemma 3.2, $g\in A=C_G(A)$ and so $A\subseteq C_G(g).$ Notice that for any element $x\in G$, we have $x\in \langle x \rangle$, which must be contained in a maximal abelian subgroup of $G$. It follows that $G\subseteq C_G(g)$ and so $g\in Z(G),$ as wanted. $\Box$ \\

For Theorem 1.4, we begin with several lemmas.\\

\noindent{\bf Lemma 3.3.} Let $G$ be a group. Suppose that $H$ satisfies the subnormalizer condition in $G$ and $M$ is a subgroup of $G.$ Then \\
(1) If $H\leq M$, then $H$ satisfies the subnormalizer condition in $M$.\\
(2) If $H$ is subnormal in $G$, then $H$ is normal in $G$.\\

\noindent{\bf Proof.} (1)  Assume that $K$ is a subgroup of $M$ such that $H$ is normal in $K.$ Since $H$ satisfies the subnormalizer condition in $G$, we have $N_G(K)\leq N_G(H)$. Thus, $M\cap N_G(K)\leq M\cap N_G(H)$ and so $N_M(K)\leq N_M(H)$, as wanted.\\

(2) Assume that $H\unlhd\unlhd G$. By induction on the length of a subnormal series from $H$ to $G$, we may suppose that  $H\unlhd K\unlhd G.$ Since $H$ satisfies the subnormalizer condition in $G$, we have $N_G(K)\leq N_G(H).$ Notice that $G= N_G(K).$ It follows that $G= N_G(H)$ and so $H\unlhd G,$ as wanted.\\

\noindent{\bf Lemma 3.4.} (1) If $H$ is a pronormal subgroup of a group $G$, then $H$ is a weakly normal subgroup of $G$.\\
(2) If $H$ is a $\mathscr{H}$-subgroup of a group $G$, then $H$ is a weakly normal subgroup of $G$.\\
(3) If $H$ a weakly normal subgroup of a group $G$, then $H$ satisfies the subnormalizer condition in $G$.\\
(4) If $H$ an $NE$-subgroup of a group $G$, then $H$ satisfies the subnormalizer condition in $G$.\\

\noindent{\bf Proof.} For (1),(2) and (3), see \cite{Ballester-Bolinches}. Now we prove (4).  Assume that $K$ is a subgroup of $G$ such that $H$ is normal in $K.$ As $H$ an $NE$-subgroup of $G$, we have $H=N_G(H)\cap H^G=N_{H^G}(H).$ Let $X=N_G(K).$ Consider $H^X$, the normal closure of $H$ in $X$. It is easy to see that $H\leq H^X\leq K$ and so $H=N_{H^X}(H)$ and $H\unlhd H^X.$ It follows that $H^X\leq N_{H^X}(H)=H$ and so $H=H^X\unlhd X.$ Hence $X=N_G(K)\leq N_G(H)$, which indicates that $H$ satisfies the subnormalizer condition in $G$. $\Box$\\

\noindent{\bf Proposition 3.5.} Let $G$ be a group. Then $\mathcal{SN}(G)=\{N_G(H)| H\in\Sigma(G) \}.$\\

\noindent{\bf Proof.} Let $S$ be a self-normalizing subgroup of $G$. Then $S=N_G(S).$  Assume that $K$ is a subgroup of $G$ such that $S$ is normal in $K.$ Clearly, $N_G(K)=N_G(S)$ and so  $H\in\Sigma(G).$ It follows that $S=N_G(S)\in \{N_G(H)| H\in\Sigma(G) \}.$ \\
\indent Let $H\in \Sigma(G)$. We show that $N_G(N_G(H))=N_G(H).$ Notice that $H\unlhd N_G(H)\unlhd N_G(N_G(H))$. Then by Lemma 3.3(1) and (2), we have  $H\unlhd N_G(N_G(H))$ and so $N_G(N_G(H))\leq N_G(H).$ Now it is obvious that $N_G(N_G(H))=N_G(H)$ and so $N_G(H)\in \mathcal{SN}(G).$  Hence, $\mathcal{SN}(G)=\{N_G(H)| H\in\Sigma(G) \}.$  $\Box$ \\

\noindent{\bf Proof of Theorem 1.4.} (a) Assume first $\Omega(G)$ is one of the first four cases. Then it is clear that $\Omega(G)$ contains all Sylow subgroups of $G$. Hence, by Lemma 3.4(1), (2) and (3), we have $\mathcal{S}(G)\subseteq \Omega(G)\subseteq \Sigma(G).$  
Thus, $N_{\Sigma}(G)\leq N_{\Omega}(G)\leq Z_{\infty}(G)$ by definition. It follows from Theorem 1.1 and Proposition 3.5, $N_{\Sigma}(G)= Z_{\infty}(G)$ and so $N_{\Omega}(G)=Z_{\infty}(G).$ \\
\indent Now assume that $\Omega(G)$ is the set of all $NE$-subgroups of $G$. By Theorem 1.1 and Lemma 3.4(4), we have $Z_{\infty}(G)\leq N_{\Omega}(G).$ Let $X=N_G(P),$ where $P$ is any Sylow subgroup of $G.$ Then we have $X=N_G(X)$, which indicates that $X$ is an $NE$-subgroup of $G.$ Hence, $N_{\Omega}(G)\leq N_G(X)=X=N_G(P).$ By the arbitrariness of $P$, we have $N_{\Omega}(G)\leq Z_{\infty}(G)$ and so $Z_{\infty}(G)=N_{\Omega}(G).$\\

(b) Notice that $\mathcal{S}_p(G)\subseteq \Omega_p(G)\subseteq \Sigma_p(G).$ Hence $N_{\Sigma_p}(G)\leq N_{\mathcal{S}_p}(G)$ by definition. Now we show that $N_{\mathcal{S}_p}(G)\leq N_G(H)$ for all $H\in \Sigma_p(G).$\\
\indent Let $H$ be a $p$-subgroup of $G$ satisfying the subnormalizer condition in $G$. Then there exists a Sylow $p$-subgroup $P$ of $G$ containing $H$. Notice that $H\unlhd\unlhd P\unlhd N_G(P).$ Then by Lemma 3.3(1) and (2), we have $H\unlhd N_G(P)$ and so $N_{\mathcal{S}_p}(G)\leq N_G(P)\leq N_G(H).$ By the arbitrariness of $H$, we have $N_{\mathcal{S}_p}(G)\leq N_{\Sigma_p}(G)$ and so (b) follows.  $\Box$\\  

\noindent{\bf Remark 3.6.} Let $\Omega(G)$ is the set of all $NE$-subgroups of $G$. The statement $N_{\Omega_p}(G)=N_{\mathcal{S}_p}(G)$ is not always true. For example,  let $G=A_5$ and $H$ be a Sylow 5-subgroup of $G$. Then we have $H^G=G$ and so $N_G(H)\cap H^G=N_G(H)\cong D_{10},$ which has order 10. So $H$ is not an  $NE$-subgroup of $G.$ It follows that $\Omega_5(G)=\emptyset$ and so $N_{\mathcal{S}_5}(G)<N_G(H)<G= N_{\Omega_5}(G).$

\section{Applications}

A group $G$ is said to be a $\mathscr{T}$-group if every subnormal subgroup of $G$ is normal in $G$. In \cite{Ballester-Bolinches}, A. Ballester-Bolinches proved that $G$ is a solvable $\mathscr{T}$-group if and only if every subgroup of $G$ satisfies the subnormalizer condition in $G.$ E. Shenkman in \cite{Schenkman} proved that the norm $N(G)$ is contained in the second center of $G$. Now it is easy to see the following:\\

\noindent{\bf Corollary 4.1.} Let $G$ be solvable $\mathscr{T}$-group. Then $Z_{\infty}(G)=N(G)\leq Z_2(G),$ where $Z_2(G)$ is the second term of the upper central series of the group $G$.



\begin{thebibliography}{0}

\bibitem{Baer1} R. Baer. Der kern, eine charakteristische untergruppe. {\it Compos. Math}, 1:254–382, 1935.


\bibitem{Baer2} R. Baer. Group elements of prime power index. {\it Trans. Amer. Math. Soc.}, 75:20–47, 1953.

\bibitem{Ballester-Bolinches} A. Ballester-Bolinches and R. Esteban-Romero. On finite $\mathscr{T}$-groups. {\it J. Aust. Math. Soc.}, 75:181–191, 2003.

\bibitem{Bianchi} M. Bianchi, A. G. B. Mauri, M. Herzog, and L. Verardi. On finite solvable groups in which normality is a transitive relations. {\it J. Group Theory}, 3:147–156, 2000.

\bibitem{Huppert} B. Huppert. {\it Endliche Gruppen I}. Springer-Verlag, Berlin etc, 1967.

\bibitem{Muller} K. H. Muller. Schwachnormale untergruppen: Eine gemeinsame verallge-meinerung der normalen und normalisatorgleichen untergruppen. {\it Rend. Sem. Mat. Univ. Padova}, 36:129–157, 1966.


\bibitem{Mysovskikh} V. I. Mysovskikh. Investigation of subgroup embeddings by the computer algebra package GAP, in: Computer algebra in scientific computing—CASC’99 (Munich), Springer, Berlin, 1999. page 309–315.


\bibitem{Schenkman} E. Schenkman. On norm of a group. {\it Illinois J. Math.}, 4(1):150–152, 1960.

\bibitem{Shirong} L. Shirong. On minimal non-$PE$-groups. {\it J. Pure Appl. Algebra}, 132:149–158, 1998.


\end{thebibliography}

\end{document}